\documentclass[8pt,openany]{article} 
\usepackage{amssymb} 
\usepackage{amsmath,amscd,mathrsfs}
\usepackage{amsmath}
\makeatletter
\renewcommand*\env@matrix[1][*\c@MaxMatrixCols c]{%
  \hskip -\arraycolsep
  \let\@ifnextchar\new@ifnextchar
  \array{#1}}
\makeatother

\usepackage{amsmath}
\usepackage{esint}
\usepackage{enumitem}
\usepackage{hyperref}
\usepackage[active]{srcltx}

\newtheorem{thm}{Theorem}

 \newtheorem{cor}{Corollary}
 
\newtheorem{prop}{Proposition}

 \def\C{{ \! \rm \ I\!\!\!C}}    \def\R{{ \! \rm \ I\!R}}
   \def\H{{ \!\! \rm \ I\!H}} 
\def\Z{{ \! \rm Z\!\!Z}}    
      
   \def \square{\hbox {$\sqcup
    $\llap {$\sqcap $}}} 
  
\newcommand{\parcial}[2]{\frac{\partial#1}{\partial#2}}

\newcommand{\normv}[1]{ \vert #1 \vert }

\title{Hopf tori and standard tori}
\author{Leonardo A. Cano
Garc\'{i}a}
\date{\small{\it Mathematics department,Universidad Nacional de Colombia,\\
    Bogot\'a, Colombia.}}
\begin{document}
\maketitle
\begin{abstract}
This article provides a complete characterization of the conformal classes of product tori and standard flat tori in complex dimension 1 (real dimension 2). Utilizing basic differential geometry methods, our approach contrasts with techniques employing Hopf tori for the conformal classification of Riemann surfaces of genus 1. While the results may be familiar to experts in complex analysis and Riemann surface theory, we contend that this work offers a clear and insightful perspective on the conformal properties of these geometrically distinct appearing tori.\end{abstract}
\section{Introduction}
\label{sec:introduction}
The upper half plane $\H \subset \C$ modulo the action of $SL(2,\Z)$ is well known to parametrize the conformal equivalence classes of tori. In this article, we focus on two important families of tori: product tori and standard flat tori. Our primary goal is to calculate their specific conformal classes within the moduli space $\mathscr{M}:=\H/SL(2,\Z)$  as presented in Theorem~\ref{thm:conformal classes rotational tori} and Theorem~\ref{thm:conformal classes product tori}. As a significant byproduct of our investigation, we show how these tori can be realized as Hopf tori (see Section~\ref{sec:hopf_tori}) whose defining curve is a circumference in the two dimensional sphere $S^2$,  and we establish the crucial result that standard tori and product tori are always conformally equivalent to each other.
\\
\\
We denote by $S^1$ the unit circle of the complex plane $\C$. Throughout this article, we define a {\bf torus}  as any manifold diffeomorphic to    $S^1 \times S^1$,  and a  {\bf complex torus} as a torus equipped with a complex structure.  Now let $R,r \in \R$ be such that $R>r$. We denote by $T_{R,r}$ the torus generated by rotating the circle of radius $r$  centered at  $(R,0)$ around the $y$--axis. The standard Riemannian metric of   $\R^3$ induces a Riemannian metric on $T_{R,r}$, whose conformal class naturally defines a complex structure. We will also use  $T_{R,r}$ to refer to this complex torus and  we term these complex tori, {\bf standard tori}.  Finally, we define {\bf product tori} as the complex tori associated to the complex structures of $S^1 \times S^1$ induced by Riemannian metrics that, in the coordinates $(\theta,\varphi) \mapsto (e^{i\theta}, e^{i\varphi})$, take the form $g=ad\theta^2+bd\varphi$ for $a,b$ real positive numbers.
\\
\\
Our primary motivation for the calculations presented in this article was to {\it explicitly visualize} the conformal classes of metrics on $S^1 \times S^1$ by finding embeddings in $\R^3$ whose induced metric is conformal to a given elliptic curve (equivalently, to a complex structure). This is directly related to a conjecture, now a theorem proven by Garsia in~\cite{Garsia} and reproved by Pinkall in~\cite{Pinkall}, which asserts that every elliptic curve can be realized via an embedding of $S^1 \times S^1$  (see \cite[Corollary in page 1]{Pinkall}). Product tori and standard tori represent natural and important families of tori whose conformal classes, as subsets of the moduli space $\mathscr{M}$, warrant a clear understanding. This article summarizes the author's journey in solving this problem. Our calculations highlight essential features of the theory of moduli spaces of conformal classes of metrics, which we believe motivates and justifies this work.
\\
\\
In this article, we demonstrate that the conformal classes of standard tori within the moduli space $\mathscr{M}$ correspond to the image of $\{0\} \times (0,\infty) \subset \H$ under the quotient map $q: \H \to \mathscr{M}$ (see Theorem~\ref{thm:conformal classes rotational tori}). As anticipated, standard tori $T_{R,r}$ sharing the same ratio $R/r$ are conformally equivalent. Specifically, we show that for every $\omega>0$  there exists a standard torus $T_{R,r}$ biholomorphic to the elliptic curve $\C/(\Z+i\omega\Z)$; in fact, the standard  torus $T_{R,r}$ corresponds to the real number $\omega:=\frac{2\pi}{\sqrt{\frac{R^2}{r^2}-1}}$. Furthermore, we establish that the conformal classes of product tori coincide with those of standard tori, thus occupying the same ray within the moduli space $\mathscr{M}$.  Moreover, we prove that these conformal classes of standard and product tori are precisely those of Hopf tori associated with circumferences in $S^2$ (see Proposition~\ref{prop:Product and Hopf}). 
\\
\\
In Dall'Acqua et al.~\cite{Dallaqua}, the authors calculate the conformal classes of what they term "tori of revolution." While closely related, these tori of revolution are distinct from our standard tori. However, Dall'acqua et al.'s~\cite[Proposition 4.2]{Dallaqua} and~\cite[Theorem 4.5]{Dallaqua} exhibit a similar spirit to our Theorems~\ref{thm:conformal classes rotational tori} and~\ref{thm:conformal classes product tori}. Indeed, as a consequence of Dall'acqua et al.'s~\cite[Theorem 4.5]{Dallaqua} and our results,  the conformal classes of tori of revolution and standard (or product) tori are the same. In Pinkall~\cite{Pinkall}, the author employed Hopf tori to construct examples of Willmore surfaces, which are extremal surfaces of the Willmore functional $\int\int H^2dA$ (where $H$ denotes the mean curvature). Also around the Willmore functional, Dallaqua et al.~\cite{Dallaqua} investigate the Willmore flow with initial data given by a torus of revolution.
\\
\\
Here we briefly recall the basic facts from the theory of elliptic curves that will be used throughout this article.  A complex torus  (also known as an elliptic curve) can be represented as the quotient of the complex plane $\C$ by a lattice $\Lambda=\omega_1\Z+\omega_2\Z$, where $\omega_1$ and $\omega_2$   are two complex numbers that are linearly independent over the reals. By scaling and rotating the lattice, we can always choose a basis $\{\omega_1,\omega_2\}$ such that $\omega_1=1$ and $\omega_2=\tau$ lies in the upper half plane $\H = \{z \in \C : \text{Im}(z) > 0\}$. The complex number $\tau$ is called the {\bf period ratio} of the complex torus (with respect to the chosen basis).
\\
\\
The moduli space of elliptic curves $\mathscr{M}$ is a central object of study in numerous areas of mathematics, spanning Riemannian geometry, complex geometry, number theory, and theoretical physics. Its significance lies, in part, in providing a unique lens through which the intricate relationships between these different fields can be observed. It is our hope that the calculations and the manner of presentation in this article will ultimately motivate or reinforce future explorations of these profound interactions. 
\\
\\
We believe that the specific calculations detailing the relationships between the conformal classes of standard tori, product tori, and Hopf tori presented in this work have not appeared in the published literature before in such a synthetic manner. Moreover, we contend that this article presents known results in a fresh and insightful way.
\\
\\
In Section~\ref{sec:calculating tau rotational}, we explicitly calculate the points in the moduli space $\mathscr{M}$ corresponding to the complex tori $T_{R,r}$. Specifically, we determine a parameter $\tau \in \H$, expressed in terms of $R$ and $r$, such that the complex torus $\C/(\Z + \tau\Z)$ is biholomorphic to $T_{R,r}$. As previously mentioned, these $\tau$ values lie along the positive imaginary axis $\{iy: y \in (0,\infty)\}$ within the upper half plane $\H$. In Section~\ref{sec:hopf_tori}, we summarize the necessary results from~\cite{Pinkall}, and we identify which Hopf tori are biholomorphic to rotational tori (see Corollary~\ref{cor: Standard and Hopf}). Finally, Section~\ref{Sec: conclusion and Pers} outlines perspectives for future work stemming from this study.
\section{The conformal class of a standard torus}
\label{sec:calculating tau rotational}
A universal covering map of the standard torus $T_{R,r}$  is given by 
\begin{equation}\label{eq: universal covering} \Phi(\theta, \varphi):=((R+r\cos(\varphi))\cos (\theta),(R+r\cos(\varphi))\sin (\theta),r\sin\varphi),
\end{equation}
where  $(\theta,\varphi) \in \R^2$. The induced Riemannian metric $g$ on $T_{R,r}$ from the standard metric of $\R^3$, when pulled back via $\Phi$ to the $(\theta,\varphi)$ coordinates induced by the covering map (\ref{eq: universal covering}), is given by
$$
g = \Phi^* (dx^2 + dy^2 + dz^2) = (R+r\cos(\varphi))^2 d\theta^2 + r^2 d\varphi^2.
$$
In these $(\theta,\varphi)$ coordinates, the metric $g$ is diagonal, indicating that $\parcial{}{\theta}$ and $\parcial{}{\varphi}$ are $g$--orthogonal. The almost complex structure $J_g$ compatible with $g$ and the orientation acts as:
\begin{equation}\label{eq: J in coord}
J_g(\partial_\theta) = \frac{R+r\cos (\varphi)}{r}\partial_\varphi \quad \text{and} \quad J_g(\partial_\varphi) = -\frac{r}{R+r\cos (\varphi)}\partial_\theta.
\end{equation}
The different coefficients of $\partial_\varphi$ and $\partial_\theta$ in the action of $J_g$ show that the coordinates $\theta,\varphi$ are not analytic in the induced complex structure. However, consider the vector fields  
$$X:=\parcial{}{\theta} \text{ and } Y:=\frac{R+r\cos (\varphi)}{r}\parcial{}{\varphi},$$
These vector fields satisfy the crucial properties 
\begin{equation}\label{eq: good orthog vect}
[X,Y]=0, \text{ } JX=Y \text{ and } JY=-X,  
\end{equation}
The vanishing Lie bracket $[X,Y]=0$ ensures, by Frobenius' theorem, the existence of local coordinates $u,v$ such that $\parcial{}{u}=X$ and $\parcial{}{v}=Y$.  In this $(u,v)$ coordinate system, the action of the almost complex structure $J_g$
becomes 
$$
J_g \parcial{}{u}=\parcial{}{v} \text{ and } J_g \parcial{}{v}=-\parcial{}{u}, 
$$
which corresponds to multiplication by $i$ on the complex coordinate $z=u+iv$. Therefore, the flows of the commuting vector fields $X$ and $Y$ define a system of analytic (holomorphic) coordinates for the complex structure induced by $g$ on the standard torus. Next we construct explicitly these coordinates.
\\
\\
Consider the function
\begin{equation}\label{eq:def f}
F(\varphi):=\int_0^\varphi \frac{r}{R+r\cos(x)}dx.
\end{equation}
Note that  the integrand includes the factor $r$  from the definition of  the vector field $Y$. The derivative of the function $F$ with respect to $\varphi$ is $F'(\varphi)=\frac{r}{R+r\cos(\varphi)}$. Since $R>r>0$, we have $F'(\varphi)>0$. This implies that $F$ is a strictly increasing function and thus has a well-defined inverse $F^{-1}$.
\\
\\
The flow of the vector field $Y=\frac{R+r\cos (\varphi)}{r}\parcial{}{\varphi}$ starting at $(\theta_0, \varphi_0)$ is given by 
$$\gamma_{(\theta_0,\varphi_0)}(t)=(\theta_0,\psi(t)),$$ where $\psi (t)$ satisfies $\psi'(t)=\frac{R+r\cos (\psi(t))}{r}$ and $\psi(0)=\varphi_0$. Notice that 
$$\frac{d}{dt} F(\psi(t))=F'(\psi (t)) \psi'(t)=\frac{R+r\cos (\psi(t))}{r} \frac{r}{R+r\cos (\varphi)}=1.$$
Integrating with respect to $t$ from $0$ to $t$, we get $F(\psi(t))=t+F(\varphi_0)$ and hence $\psi(t)=F^{-1}(t+\varphi_0)$.
\\
\\
To construct explicit analytic coordinates, let us consider the map  $\Psi:\R^2 \to \R^2$ defined by taking the flows with initial point $(0,0)$. Let  $$\Psi(\theta,s):= (\theta,F^{-1}(s)).$$
This map is a global chart of the plane $(\theta, \varphi)$.  We claim that $z=\theta+it$ are analytic coordinates for the complex structure induced by the almost complex structure $J_g$. Indeed, we have
$$
\parcial{}{\theta} \Psi= \parcial{}{\theta}= X \text{ and }\parcial{}{t} \Psi=\parcial{}{t}(\psi) \parcial{}{\varphi}=Y.
$$
Thus, in the $\theta, t$ coordinates,  the vector fields $\parcial{}{\theta}$ and $\parcial{}{t}$  correspond to $X$ and $Y$ respectively, and we know that $J_gX=Y=\partial_t$ and $J_gY=-X=-\partial_\theta$. This confirms that $J_g \partial_\theta=\partial_t$ and $J_g \partial_t=-\partial_\theta$, which is the standard action of the complex structure on the real and imaginary parts of the analytic coordinate $z=\theta+it$. 
\begin{prop}\label{prop:Pi covering}
The map $\Pi: (\theta, s) \mapsto \Phi(\Psi(\theta,s))$ is a holomorphic universal covering map from the complex plane $\C \cong \R^2$   with the complex structure $z=\theta+is$ to the complex torus $T_{R,r}$ with the complex structure induced by its Riemannian metric. 
\end{prop}
{\bf Proof:}
\\
The map  $\Phi:\R^2 \to T_{R,r}$ is  smooth and surjective, locally a diffeomorphism, and its fibers are discrete, so it is easy to check that  it is a universal covering map. The function $\Psi$ is a diffeomorphism  from $\R^2$ to $\R^2$, thus, their composition $\Pi:=\Phi \circ \Psi$ is also a smooth surjective map, locally a diffeomorphism, and a universal covering map. 
\\
\\   
Let $J_{\C}$ be the standard complex structure on $\C$ such that $J_{\C}(\partial_\theta) = \partial_s$ and $J_{\C}(\partial_s) = -\partial_\theta$. Let $J_T$  be the complex structure on $T_{R,r}$ induced by its Riemannian metric $g$. The map $\Pi$ is holomorphic if and only if $d\Pi \circ J_{\C} = J_{T} \circ d\Pi$, which we will check next.
\\
\\
We know that $\Psi(\theta,s)=(\theta, F^{-1}(s))$, so the differential of $\Psi$ maps $\parcial{}{\theta}$ to $\parcial{}{\theta}$ and $\parcial{}{s}$ to the vector field $Y$ defined above. The differential of  $\Phi$ maps $\parcial{}{\theta}$ and $\parcial{}{s}$ to vector fields of $T_{R,r}$.  Let $Z_1:=\Phi_*(\parcial{}{\theta})$ and $Z_2:=\Phi_*(\parcial{}{\varphi})$ be a local frame on $T_{R,r}$. One can check from (\ref{eq: J in coord}) that 
\begin{equation}\label{eq: JT en Z1 y Z2}
J_T(Z_1)=\frac{R+r\cos (\varphi)}{r}Z_2 \text{ and }J_T(Z_2)=-\frac{r}{R+r\cos (\varphi)}Z_1,
\end{equation}
Now let us  examine
\begin{equation}\label{eq: dPi en Z1 y Z2}
d\Pi(\parcial{}{\theta})=d \Phi(d\Psi(\parcial{}{\theta}))=Z_1 \text{ and } d\Pi(\parcial{}{s})=d \Phi(d\Psi(\parcial{}{s}))=\frac{R+r\cos(\varphi)}{r}Z_2.
\end{equation}
One can apply equations (\ref{eq: JT en Z1 y Z2}) and (\ref{eq: dPi en Z1 y Z2}) to prove $d\Pi \circ J_{\C} = J_{T} \circ d\Pi$. $\square$
\\
\\
Let $a:=R/r$, and define
$$
\omega := \omega(a) := \int_0^{2\pi} \frac{r}{R+r\cos(x)}dx = \int_0^{2\pi} \frac{1}{a+\cos(x)}dx = \frac{2\pi}{\sqrt{a^2-1}}.
$$
By the definition of the function $F$ (see (\ref{eq:def f})), we have the property that
$$
F(\varphi+2m\pi) = F(\varphi) + m\omega,
$$
for every $m \in \Z$. This periodicity of F leads to a corresponding periodicity of the map $\Pi:=\Phi \circ \varphi$. For all $m,n \in \Z$, we have
\begin{equation}\label{eq: def P}
\Pi(\theta+2\pi m, s+n\omega)=\Phi(\Psi(\theta+2\pi m, s+n\omega)) = \Phi(\theta+2\pi m, F^{-1}(s+n\omega)).
\end{equation}
Since $\Phi$ has period $2\pi$ in $\theta$ variable, and $F^{-1}(s+n \omega)=F^{-1}(F(\varphi)+n \omega)=\varphi+2 n \pi$ where $s=F(\varphi)-F(0)=F(\varphi)$,  we have $$
\Pi(\theta+2\pi m, \varphi + 2n\pi) = \Pi(\theta, \varphi).
$$
Therefore,
\begin{equation}\label{eq: def P revised}
\Pi(\theta+2\pi m, s+n\omega) = \Pi(\theta, s).
\end{equation}
Then, (\ref{eq: def P revised}) implies that we can define a map $H:\C/(2\pi \Z+i \omega  \Z) \to T_{R,r}$ by
$$
H([\theta+is]) := \Pi(\theta,s),
$$
where $[\theta+is]$ denotes the equivalence class of the complex number $\theta+is$ in the quotient $\C/\Lambda$ with lattice $\Lambda:=2\pi \Z+i\omega  \Z$. From (\ref{eq: def P revised}) and the fact that $\Pi$ is a holomorphic universal covering map (see Proposition~\ref{prop:Pi covering}), it follows that $H$ is a well-defined biholomorphic map between the complex torus $\C/(2\pi \Z+i\omega  \Z)$ and the standard torus $T_{R,r}$.
\\
\\
The following theorem summarizes the results proved with the previous calculations.
\begin{thm}\label{thm:conformal classes rotational tori}
The standard torus $T_{R,r}$ is biholomorphic to the complex torus $\C/(2\pi \Z+i\omega(R/r) \Z)$, where $\omega(R/r)=\frac{2\pi}{\sqrt{(\frac{R}{r})^2-1}}$. 
\end{thm}
An easy and expected corollary of the previous theorem is the following.
\begin{cor}
if $R/r=R'/r'$ then the two standard tori $T_{R,r}$ and $T_{R',r}$ are conformally equivalent.
\end{cor}
The range of the function $\omega(a):=\frac{2\pi}{\sqrt{\frac{R^2}{r^2}-1}}$ for $a= R/r>0$ is the open interval $(0,\infty)$. This implies that standard tori are biholomorphic to elliptic curves associated with period ratios $\tau=iy $ where $y>0$, which corresponds to the positive imaginary axis $\mathscr{R}:=\{iy:y>0\}$ in the upper half plane $\H$. However, the moduli space of elliptic curves is given by the quotient $\mathscr{M}:=\H/SL(2,\Z)$  and a fundamental domain for this action (see  \cite[Lemma 2.3.1]{DiamondShurm}) is the set
$$
D := \{x+iy \in \H : -1/2 \leq x \leq 1/2 \text{ and } |x+iy| \geq 1 \}.
$$
The existence of this fundamental domain implies that different points in $\H$, and therefore potentially different values of $\omega$, can represent the same conformal class of torus. We will now describe these additional conformal equivalences between standard tori arising from the $SL(2,\Z)$ action on $\H$.
\\
\\
The conformal equivalence between the elliptic curves $\C/(\Z+\tau\Z)$ and $\C/(\Z-\frac{1}{\tau}\Z)$ arises from the fact that the M\"obius transformation $\tau \mapsto \frac{-1}{\tau}$ is associated with the matrix $\begin{pmatrix}
0&-1\\
1&0
\end{pmatrix}$ which belongs to $SL(2,\Z)$. This transformation establishes a conformal equivalence, and in particular, it provides a bijective correspondence between the conformal classes of elliptic curves associated with period ratios in the regions
$\{\tau \in \H: {\rm Re}(\tau)<1\}$ and $\{\tau \in \H: {\rm Re}(\tau)>1\}$. The following proposition expresses this conformal equivalence in terms of standard tori.
\begin{prop}
If  $1<R/r<\sqrt{2}$, then the standard torus $T_{R,r}$ is conformally equivalent to the standard torus $T_{R,\sqrt{R^2-r^2}}$. Moreover, the map $T_{R,r} \mapsto T_{R,\sqrt{R^2-r^2}}$ is a bijection between the set of standard tori  $T_{R,r}$ such that $1<R/r<\sqrt{2}$,  and the set of standard tori  $T_{R,r}$ such that $\sqrt{2}<R/r$.
\end{prop}
{\bf Proof}:
Recall that $a:=R/r$.  The imaginary part of the period ratio of the complex torus biholomorphic to $T_{R,r}$ is $\omega_a=\frac{1}{\sqrt{a^2-1}}$.  The inequality  $0<\omega_a \leq 1$ is equivalent to $\sqrt{2}<a$. Thus, standard tori with
$\sqrt{2}<R/r$ are biholomorphic to complex tori with imaginary part of the period ratio in $(0,1]$. Similarly, the inequality  $\omega_a \geq 1$  is equivalent to $1<a<\sqrt{2}$. Thus, standard tori with $1<R/r< 
2$ are biholomorphic to complex tori with imaginary part of the period ratio in $(1,\infty)$.
\\
\\
Now, consider the torus $T_{R,r'}$ where $r':=\sqrt{R^2-r^2}$. Let $b=R/r'=\frac{a}{\sqrt{a^2-1}}$. So the imaginary part  of the period ratio for $T_{R,r'}$ is $\omega_b=\frac{1}{\omega_a}$. Because of this, the complex tori $\C/(2\pi \Z + i\omega_a \Z)$ and $\C/(2\pi \Z + i\omega_b \Z)$ are conformally equivalent because their period ratios $\tau_a:=\frac{\omega_a}{2\pi}$ and $\tau_b:=\frac{\omega_b}{2\pi}$ are related by the  $SL(2,\Z)$ M\"obius transformation $\tau \mapsto \frac{-1}{\tau}$. Finally, the bijection follows from the function  $f(a):=\frac{a}{\sqrt{a^2-1}}$.  For $1<a< 
\sqrt{2}$, we have $f(a)>\sqrt{2}$, and this function provides a bijection between the intervals $(1, 
\sqrt{2})$ and $(\sqrt{2},\infty)$ as its derivative is negative, indicating it is strictly decreasing.  Therefore, the map $T_{R,r} \mapsto T_{R,\sqrt{R^2-r^2}}$  induces a bijection between the specified sets of standard tori. $\square$
\section{The conformal classes of product tori}\label{Sec:product metric}
For $0<a\leq 1$, let $S_a$ be the product torus $S^1\times S^1$ endowed with the Riemannian metric which in coordinates $(\theta,\varphi) \mapsto (e^{i\theta}, e^{i\varphi})$ takes the form
$$
g_a=d\theta^2+a d\varphi^2.
$$
Every product torus (with a metric of the form $bd\theta^2+cd\varphi^2$) where $b,c>0$) is conformally equivalent to a torus of the form $S_a$ for some $a \in [1,\infty)$. This can be seen by scaling one of the coordinates by $\sqrt{c/b}$ (or $\sqrt{c/b}$).
\\
\\
Recall that the metric $g_a$ induces a complex structure $J_a$ on $S_a$. The vector fields $X=\partial_\theta$ and $\frac{1}{a}\partial_\varphi$ are orthogonal and satisfy $$[X,Y]=0 \text{ and }  g(X,X)=g(Y,Y).$$ Furthermore, $X$ and $Y$ are $g_a$--orthonormal.  Proceeding as in Section~\ref{sec:calculating tau rotational}, their flows define analytic coordinates. The map $$G([\theta+is])=(e^{i \theta},e^{i\frac{1}{a}s})$$
defines a biholomorphism from $\C/(\Z+ \frac{i}{a}\Z)$ to $(S^1 \times S^1,J_a)$.  These calculations prove the following theorem.
\begin{thm}\label{thm:conformal classes product tori}
Let $a>b>0$. Let $S_{a,b}$ be the product torus $S^1\times S^1$ endowed with the Riemannian metric which in coordinates $(\theta,\varphi) \mapsto (e^{i\theta}, e^{i\varphi})$ takes the form
$$
g_a=bd\theta^2+a d\varphi^2.
$$
Then $S_{a,b}$ is conformally equivalent to $\C/(\Z+ \frac{bi}{a}\Z)$.
\end{thm}
As a consequence of Theorems~\ref{thm:conformal classes rotational tori} and~\ref{thm:conformal classes product tori} we have
\begin{cor}\label{Cor:Standard and prod tori}
The standard torus $T_{R,r}$ is conformally equivalent to  the product torus $S_a$ where $a=\sqrt{\frac{R^2}{r^2}-1}$. 
\end{cor} 
\section{Hopf Tori}
\label{sec:hopf_tori}
First we summarize some  relevant results of~\cite{Pinkall}. 
Let $Q$ denote the quaternions and let $S^3$ be all the quaternions of norm $1$. Let us indentify $S^2$ with the unit sphere in the real subspace of $Q$ spanned by $1,j$ and $k$.  Let $q \in Q$ and suppose that $q=a+ib+cj+dk$, for $a,b,c,d \in \R$, then
$$
\tilde{q}:=a-ib+cj+dk
$$
is an antiautomorphism of $Q$. The version of the Hopf fibration $\pi:S^3 \to S^2$ used in~\cite{Pinkall} is defined by $$\pi(q):=\tilde{q}q.$$ With the $S^1$--action on the fibers defined by $$(e^{i\varphi},q) \mapsto e^{i \varphi}q,$$ the  Hopf fibration $\pi:S^3 \to S^2$ is a $S^1$--principal bundle. 
\\
\\
By definition a {\bf Hopf torus} is the inverse under the projection $\pi$ of a closed curve of $S^2$.  Let $p:[a,b] \to S^2$ be a smooth Jordan curve on $S^2$. We define the associated Hopf torus as $$T_p:=\pi^{-1}(p([a,b])),$$ 
We consider $T_p$  as a complex torus endowed with the conformal class of the Riemannian metric obtained by pulling back the standard Riemannian metric of $S^3$ via the inclusion $T_p \hookrightarrow S^3$.
\\
\\
Let $\alpha$ be the connection $1$--form of of the principal $S^1$--bundle $\pi:S^3 \to S^2$, associated with the Riemannian metric  $g$ of $S^3$ inherited from  the quaternion norm. The horizontal bundle $H:=V^\perp$ is the orthogonal complement of the vertical bundle $V:={\rm Ker}(\pi)$ with respect to $g$.  Let $\eta:[0,L] \to S^3$ be any $\alpha$--horizontal lift of   the Jordan curve  $p$, parametrized by arc length $t$ (where $L$ is the length of $\eta$). The function
$$
\chi(t + i\varphi) = e^{i \varphi} \eta(t)
$$
is a universal covering map of the Hopf torus $T_p$. Moreover, $\chi$  is a holomorphic function from $\C$ to $T_p$ (see~\cite[Proposition 1]{Pinkall}) because the induced metric on $\C$ via $\chi$ satisfies 
\begin{equation}\label{eq: Cauchy Riemann for chi}
g(\partial_t \chi,\partial_t \chi)=1=g(\partial_\varphi \chi,\partial_\varphi \chi)=1 \text{ and } g(\partial_\varphi \chi,\partial_t \chi)=0.
\end{equation}
These equalities show that the induced metric is conformal to the Euclidean metric on $\C$, and thus $\chi$ is holomorphic.
\\
\\
Since $\eta$ is parametrized by arc length, its tangent vector $\eta'$ has unit norm and is orthogonal to $\eta$ (when considering $\eta$ as a vector in $Q \cong R^4$). Following Pinkall~\cite{Pinkall},  there exists a function $u:[a,b] \to {\rm span}\{j,k\}$ with $\normv{u}=1$ such that $$\eta'=u\eta.$$
The curve $\pi \circ \eta$ is a reparametrization of $p$. Abusing notation, we will continue to denote $\pi \circ \eta$ by $p$. As pointed out in~\cite[Formulas (6)]{Pinkall}, the tangent vector to $p$ is given by
$$
p'(t) = 2 \tilde{\eta}(t) u(t) i \eta(t).
$$
This implies that if $L$ is the length of the Jordan curve p, then $L/2$ is the length of the curve $\eta$.
\\
\\ 
Since $\eta(L/2)$ and $\eta(0)$ both lie in the fiber $\pi^{-1}(p(0))$, there exists a real number $\delta$ such that
$$
\eta(L_p/2) = e^{i \delta} \eta(0).
$$
Pinkall proves in~\cite[Proposition 1]{Pinkall} that 
$$\delta=A(p)/2,$$
where $A(p)$ is the (signed) enclosed area of the curve $p$, which we now define. Let $dV$ be the canonical volume form of $S^2$.  The {\bf signed enclosed area of $p$} is then naturally defined as
\begin{equation}\label{eq:enclosed area}
A(p) = \int \int_c  dV,
\end{equation}
where $c$ is chosen such that $\int \int_c dV \in [-2 \pi, 2\pi)$. 
\\
\\
The Gaussian curvature of $S^2$ is $1$ at every point. If $c$ is one of the two connected components of $S^2-p$, then by the Gauss-Bonnet theorem (applied to the disk $c$), we have:
$$
\text{Vol}(c) = \int \int_c dV = 2\pi - \int_p k_gds
$$
where $k_g$ denotes the geodesic curvature of $p$, and the integral is taken over one orientation of $p$.
The signed enclosed area of $p$ can also be expressed in terms of the geodesic curvature as
\begin{equation}\label{eq:enclosed area}
A(p)=\begin{cases}
2\pi-\int_{p}k_gds&\text{ if }2\pi-\int_{p}k_gds<2 \pi.\\
2\pi+\int_{p}k_gds&\text{ otherwise.}
\end{cases}
\end{equation}
Building upon the summarized results from Pinkall~\cite{Pinkall}, we now present the following theorem:
\begin{thm}\label{thm:classification Hopf tori}{(\upshape c.f \cite[Proposition 1]{Pinkall})}
The complex torus $T_p$ associated to a Jordan curve $p:[a,b] \to S^2$ is biholomorphic to $\C/(2\pi  \Z+(\frac{A}{2}+i\frac{L}{2}) \Z)$ where $L$ is the lenght of the curve $p$ and $A$ is the area enclosed by $p$.
\end{thm}
According to~\cite[Equation (15)]{Pinkall}, the isoperimetric inequality on $S^2$ implies that the signed enclosed area $A(p)$ and the length $L(p)$ of any Jordan curve $p$ satisfy
\begin{equation}\label{eq:isoperimetric ineq}
(A(p)/2 - \pi)^2 + (L(p)/2)^2 \geq \pi^2.
\end{equation}
In the plane with coordinates $A/2$, $L/2$ , the region defined by (\ref{eq:isoperimetric ineq}) includes the set
$$\{a+ib:b>0, 0\leq a\leq \pi\},$$
 This region, when combined with its reflection across the $a=0$ axis, covers a fundamental domain of the action of $SL(2,\Z)$  on the upper half plane $\H$. Since the orientation of a closed curve on $S^2$  changes the sign of the enclosed area $A(p)$ (as per the definition), Hopf tori exhaust all possible conformal classes of Riemann surfaces of genus $1$.
\\
\\
Via stereographic projections, every Hopf torus can be mapped conformally to a torus embedded in $\R^3$. Therefore, Theorem~\ref{thm:classification Hopf tori} implies the following theorem, which relates to Garsia's Conjecture: 
\begin{thm}\label{thm:GarsiaConjecture}{\upshape (c.f \cite[Theorem 1]{Pinkall})}
Let $X$ be a  compact Riemann surface of genus $1$ and let us denote $g_0$ the Riemannian metric of $\R^3$. Then, there exists an embedding $\varphi: S^1 \times S^1 \to \R^3$ such that $X$ and $(S^1 \times S^1, \varphi^*g_0)$ are conformally equivalent.
\end{thm}   
Pinkall proves in~\cite{Pinkall} that the embedding whose existence is claimed in Theorem~\ref{thm:GarsiaConjecture} can be chosen as an algebraic surface of degree eight (see \cite[Corollary, page 381]{Pinkall}). In Sections~\ref{Sec: hopf and product} and \ref{Sec: hopf and rotational} we will determine the curves on $S^2$ associated with product tori and standard rotational tori via the Hopf fibration.
\subsection{Hopf tori and product tori}\label{Sec: hopf and product}
For $0 <t\leq 1/2$, let $H_t$  denote  to the product torus $S^1 \times S^1$ equipped with the Riemannian metric that, in the coordinates $(\theta,\varphi) \mapsto (e^{i\theta}, e^{i\varphi})$  takes the form
$$
g_t = t^2 d\theta^2 + (1-t^2) d\varphi^2.
$$
Every flat product torus (with a metric of the form $Ad\theta^2 + Bd\varphi^2$ where $A,B$ are real positive numbers) is conformally equivalent to a torus of the form $H_t$ for some $t \in  (0,1/2]$. This can be seen by scaling the coordinates. The product torus $H_t$ can be isometrically embedded in the quaternions $Q$ by
$$
\Phi(e^{i \theta}, e^{i \varphi}) = t \cos(\theta) + t \sin(\theta) i + \sqrt{1-t^2} \cos(\varphi) j + \sqrt{1-t^2} \sin(\varphi) k.
$$
With a slight abuse of notation, we will also denote the image of $\Phi$ by $H_t$. This torus $H_t$ is invariant under the diagonal $S^1$ action:
$$
e^{is} \Phi(e^{i \theta}, e^{i \varphi}) = \Phi(e^{i (\theta+s)}, e^{i (\varphi+s)}).
$$
Since we have
$$\pi(\Phi(\theta,\varphi))=\tilde{\Phi}\Phi(\theta, \varphi)=2t^2-1+2t\sqrt{1-t^2}\left(\cos(\theta-\varphi)j+\sin(\theta-\varphi)k\right),$$
the image of $S_t$ under the Hopf projection $\pi(q)= \tilde{q}q$ is the circle $p_t$  on $S^2$ parametrized by 
\begin{equation}\label{eq: circle Prod Hopf}
s\mapsto 2t^2-1+2t\sqrt{1-t^2}\left(\cos(s)j+\sin(s)\right).
\end{equation}
It is straightforward to see that $S_t$ is the Hopf torus associated with this circle $p_t$. Therefore, we can use Theorem~\ref{thm:classification Hopf tori} to calculate its conformal class as a point in the moduli space $\mathscr{M}$.
\\
\\
The signed area enclosed by $p_t$ is 
$$
A(t)=4\pi t^2
$$ 
and its lenght is 
$$
L(t)=4\pi t \sqrt{1-t^2},
$$
Therefore, by Theorem~\ref{thm:classification Hopf tori}, the complex torus $S_t$ is biholomorphic to
$$
\C / (2\pi \Z + (\frac{4\pi t^2}{2} + i \frac{4\pi t \sqrt{1-t^2}}{2}) \Z) = \C / (2\pi \Z + (2\pi t^2 + 2\pi i t \sqrt{1-t^2}) \Z).
$$
Dividing the lattice by $2\pi$, we get the conformally equivalent torus $\C / (\Z + (t^2 + i t \sqrt{1-t^2}) \Z)$.  Applying the M\"obius transformation $\tau \mapsto -1/\tau$ , we obtain
$$
-\frac{1}{t^2 + i t \sqrt{1-t^2}} = -\frac{t^2 - i t \sqrt{1-t^2}}{t^4 + t^2(1-t^2)} = -\frac{t^2 - i t \sqrt{1-t^2}}{t^2} = -1 + i \frac{\sqrt{1-t^2}}{t}.
$$
This period ratio  $-1 + i \frac{\sqrt{1-t^2}}{t}$ is $SL(2,\Z)$--equivalent to $i \frac{\sqrt{1-t^2}}{t}$ by adding 1. Hence, the product torus $S_t$  is conformally equivalent to $\C / (\Z + i \frac{\sqrt{1-t^2}}{t} \Z)$, which matches the result claimed in Theorem~\ref{thm:conformal classes product tori} (with the identification $t=b/a$ or $t=a/b$ depending on the parametrization). We have thus shown the following proposition.
\begin{prop}\label{prop:Product and Hopf}
The conformal classes of flat product tori coincide with the conformal classes of Hopf tori associated with circles on $S^ 
2$. Specifically, given $a>1$,  the flat product torus  $S_a$ (with metric $d\theta^2+ad\varphi^2$)  is conformally equivalent to the Hopf torus  $H_{t}$ associated to the circle parametrized by (\ref{eq: circle Prod Hopf}) where $t=\frac{1}{\sqrt{a^2+1}}$.
\end{prop}
\subsection{Hopf tori and standard tori}\label{Sec: hopf and rotational}
Let $R>r>0$ and suppose that $R/r>\sqrt{2}$ and define $b:=\sqrt{\frac{R^2}{r^2}-1}>1$. Corollary~\ref{Cor:Standard and prod tori} implies that the standard torus $T_{R,r}$ is conformally equivalent to the product torus $S_{b}$ (with metric $d\theta^2+bd\varphi^2$),  and Proposition~\ref{prop:Product and Hopf} implies that $S_{b}$ is conformally equivalent to $H_{t(b)}$ where $t=\frac{1}{\sqrt{b^2+1}}=\frac{r}{R}$.  We have thus proved the following corollary of the previous sections.\begin{cor}\label{cor: Standard and Hopf}
Let $R>r>0$ and suppose that $R/r>\sqrt{2}$. The standard torus $T_{R,r}$ is conformally equivalent to the Hopf torus $H_t$ associated to the circle in $S^2$ parametrized by (\ref{eq: circle Prod Hopf}) where $t=\frac{r}{R}$.
\end{cor}
\section{Conclusion and perspectives}\label{Sec: conclusion and Pers}
This article determined the conformal classes of standard, product, and Hopf tori. We showed standard and product tori are mutually conformally equivalent, and further demonstrated that the conformal classes of flat product tori coincide with those of Hopf tori associated with circles on $S^2$, explicitly linking their defining parameters. This work clarifies the relationships between these geometric tori within the moduli space of genus one Riemann surfaces.
\\
\\
Looking ahead, the theory of moduli spaces is often approached with the tools of algebraic geometry, in contrast to the differential geometric perspective adopted in this article. It would be interesting to explore whether the insights and structures from the algebraic geometric viewpoint can be translated or related to the realm of Hopf tori, which possesses a more pronounced differential geometric flavor.
\section*{Aknowledgments}
A portion of this work was completed at the Max Planck Institute of Mathematics in Bonn, and I gratefully acknowledge their hospitality. I extend my sincere thanks to Matilde Marcolli for initially suggesting Pinkall's article many years ago, and to Andrés Ángel for bringing it to my attention again. I am also indebted to Sergio Carrillo for his insightful discussions on aspects of this project, and to Alexander Cruz and Andrés Villaveces for their encouragement and for inspiring me to pursue this writing.
\bibliographystyle{alpha}
\bibliography{literatur}
\end{document}